\documentclass[12pt]{article}
\usepackage{amsmath}
\newcommand{\bee}{\begin{eqnarray*}}
\newcommand{\ene}{\end{eqnarray*}}
\newcommand{\beeq}{\begin{equation}}
\newcommand{\eneq}{\end{equation}}
\newtheorem{lem}{Lemma}[section]
\newcommand{\bel}{\begin{lem}}
\newcommand{\enl}{\end{lem}}
\newtheorem{exap}{Example}[section]
\newcommand{\beex}{\begin{exap}}
\newcommand{\enex}{\end{exap}}
\newtheorem{theo}{Theorem}[section]
\newcommand{\beth}{\begin{theo}}
\newcommand{\enth}{\end{theo}}
\newtheorem{prop}{Proposition}[section]
\newcommand{\bep}{\begin{prop}}
\newcommand{\enp}{\end{prop}}
\newtheorem{cor}{Corollary}[section]
\newcommand{\bec}{\begin{cor}}
\newcommand{\enc}{\end{cor}}
\newtheorem{rem}{Remark}[section]
\newcommand{\ber}{\begin{rem}}
\newcommand{\enr}{\end{rem}}
\newtheorem{defi}{Definition}[section]
\newcommand{\bef}{\begin{defi}}
\newcommand{\enf}{\end{defi}}

\def\ci{\perp\!\!\!\perp}
\usepackage{amsmath}
\usepackage{graphicx}

\begin{document}
\title { 
SIMULTANEOUS CAUSAL INFERENCE\\ FOR\\ MULTIPLE TREATMENTS\\ VIA \\ SUFFICIENCY
}
\author {Yannis G. Yatracos \\  Cyprus U. of Technology}
\maketitle
\date{}
{\em e-mail:} yannis.yatracos@cut.ac.cy
\pagebreak
\vspace{2in}
\begin{center} \vspace{0.05in} {\large Summary} \end{center}
\vspace{0.03in}
\begin{center}
\parbox{4.8in}
{\quad  
 Some units from a population receive the same  treatment  that is different from treatments available  for other reservoir populations. The minimal sufficient statistic $s$  for the pre-treatment $x$-covariates's distributions  in  the populations
 is the coarsest balancing score. $s$ is  used to select matching units for {\em simultaneous} causal comparisons of {\em multiple} treatments.Necessary and sufficient conditions on the posterior distribution of the treatment variable (given $x$)  determine whether a statistic is either sufficient or minimal sufficient for the x-covariates' distributions. Results in the literature are thus extended. Strong ignorability of treatment assignment given $s(x)$ is also established. Consequently,  the expected treatments' differences given $s(x)$ are shown to be simultaneously unbiased for the average causal effects of {\em all}  treatments' differences. The existing statistical theory for $s$ and its estimates support their use in causal inference.
  
}
\end{center}

\vspace{.05 in}

\bigskip
{\it Some key words:} \quad    Causal Inference, Coarsest Balancing Score, Generalized  Propensity Scores, Generalized Linear Models, Matching, Minimal Sufficient 
Statistic, Propensity Function

{\it Running head:} \quad S-matching for Simultaneous Causal Inference 




\pagebreak

\section {Introduction}

\quad 
When  reservoir populations receive each a different treatment,
 the minimal sufficient statistic $s$ of the  pre-treatment $x$-covariates' distributions 
 is  used to select matching units for 
{\em simultaneous} causal comparisons of the treatments. Strong ignorability of treatment assignment 
given $s(x)$ is established and
the expected treatments' differences  given any
$s$-value  are shown to be {\em simultaneously}  unbiased for the average causal effects of
{\em all}  treatments' differences.
 Criteria are provided to obtain either  $s$ 
or a balancing score 
from the posterior distribution, $q(t|x),$ of the treatment variable, $T,$  given the $x$-covariates.
The results in Imai and van Dyk (2004) are extended, providing balancing scores for a larger 
class of $q(t|x)$-models.

 
For  two treatments,
$t=1,\  2,$
 Rosenbaum and Rubin (1983) propose 
the 
scalar propensity score $e(x)$
to balance the  pre-treatment covariates, $x  (\in R^p),$ of the
$n$ units in the treatment groups;
$e(x)$ is the conditional probability of receiving, say,  treatment 1 given $x.$
It is stated therein that $e(x)$ is the coarsest balancing score
 and it is showed, 
among others, that
if treatment assignment
and the potential units' responses 
to treatments,  $r_i(1)$ and $r_i(2), \ i=1,\ldots,n,$ 
 are conditionally independent given $x,$ 
then the difference between the 
sample treatments' means given 
$e(x)$ is
unbiased
for the average causal effect $E\{ r(2)-r(1) \}; \ E$ denotes expectation 
over the whole population.

R. Bahadur recognized that
 $e(x)$ 
is equivalent to  the likelihood ratio 
 of the $x$-populations' densities 
which is minimal sufficient (Rubin and Thomas, 1996, p. 250).
As we found recently,  the minimal sufficient statistic $s$ is used with more than two treatments 
{\em a)}  for dimension reduction theory, in particular when propensities do not exist 
(Nelson and Noorbaloochi, 2009), and  
{\em b)} in  causal inference  for each pair of treatments  {\em assuming} that 
strong ignorability of treatment assignment holds (Noorbaloochi, Nelson and Asgharian, 2010, p. 12, lines -2, -1, p. 13, lines 1-4).
This assumption and the use of  the $x$-covariates distributions  to obtain $s$  constitute two of the differences with this work.

 
For more than two treatment levels 
 Joffe and Rosenbaum (1999) study 
 causal effects using
 a small number of balancing linear functions of $x.$ 
 For multi-valued categorical treatments 
 Imbens (2000)  introduced the generalized  propensity score $P(T=t|x)$  
and
used it to estimate average causal effects for treatments' pairs  {\em but not}  for  simultaneous  causal comparisons of all
treatments; see Imbens (2000,  page 709, lines 15 to 17 and lines -4 to -1).  
For general treatment regimes 
with any type of $t$-values,
Imai and Van Dyk (2004) 
introduce the propensity function $e_{\psi}(\cdot |x),$ that has the form of $q(t|x),$  and assume that for all $x$-values it  depends on $x$ only 
through a unique, finite dimensional parameter $\theta_{\psi}(x)$ ({\em Assumption 3}, p. 856);  $\psi$ is a known parameter. 
They show  that $e_{\psi}(\cdot |x)$ (i.e. 
$\theta_{\psi}(x)$)  is  balancing score ({\em Result 1}, p. 856),  {\em not necessarily the coarsest}, and use it for simultaneous causal comparisons.
Their {\em Assumption 3}, that does not hold when $q(t|x)$ belongs  to a general exponential family,  and $\theta_{\psi}$'s dimensionality  constitute
two of the differences with this work where $s$  is not necessarily finite-dimensional and is obtained without additional assumptions.
Corollary \ref {c:extImVD}  shows also  that  {\em Assumption 3} can be weakened to allow for more $q(t|x)$ models
and  {\em Result 1}
 holds automatically.

An estimate of $s$  may not be sufficient and the same also  holds for estimates of the propensity score, the generalized propensity score and the propensity function.  This  has been neglected so far in the Causal Inference  literature
that has not its own  tools to confirm the balancing property unlike $s;$ see section 4.



The findings 
explain clearly what ``matching''  means. 
Units from different populations receiving each a different treatment form a matching group
when they 
provide the same information for the $x$-covariates' distributions.
Such groups can be used in causal comparisons, for example, 
to determine the ``right'' dose 
for a new drug, by examining {\it simultaneously} the expected response differences 
$E\{r(t_2)-r(t_1)|s(x)=s_0\}, \   E\{r(t_3)-r(t_2)| s(x)=s_0\}, \ldots, \\  E\{(r(t_k)-r(t_{k-1})|s(x)=s_0\}$
for different doses' levels $t_1 < t_2< \ldots < t_k.$ 



The framework is presented in section 2. The main theoretical results and $s$ for
generalized linear models are in section 3.
In section 4 some directions are  given for  $s$-matching's implementation in practice.
The proofs are in the Appendix.


\section{Causal inference framework and assumptions}

\quad For a random vector $U$ use $p_U(u)$ to denote its density (but also its probability). When random vector  $V$ 
is also available use $p(u|v)$ to 
denote the conditional density of $U$ given $V.$  
Let ${\cal T}$ denote the treatments and let 
$T$ be the treatment variable with values $t$
in ${\cal T}$  and prior density $\pi_T.$  Treatment 
$t$ is  used in selected units of  population
${\cal P}_t$  having  balanced $x$-covariates
with respect to ${\cal T}.$ 
The units in ${\cal P}_t$ have covariates 
$x \in {\cal C}({\cal P}_t) \subset R^p$ and
unless otherwise stated it is assumed that ${\cal C}({\cal P}_t)={\cal C}, 
t \in {\cal T}.$
Let 
$p(x|t)$ denote the $x$-covariates' density
of units in ${\cal P}_t$ and
let ${\cal D}_{\cal T}=\{p(x|t), \ t \in {\cal T}\}; \ p_X(x)$ is the marginal density of the $x$-covariates. 
The notation $p(x|t)$ does not
mean necessarily that $p$ is the same density  with the parameter $t$ changing,  $t  \in {\cal T},$
but simply denotes  the covariates' distribution in ${\cal P}_t.$
Use $q(t|x)$ to denote $T$'s density (or probability) given the $x$-covariates.
For unit $i,  
\ r_i(t)$ is
the response for treatment $t$ and
the potential outcomes 
${\cal R}$ is the set $\{r_i(t), t \in {\cal T},  \mbox{ for } i=1,\ldots, n\}.$ 
Conditional independence of $x$ and $y$ given $z$ is denoted by 
$ x \ci y| z$ (Dawid,1979).
The expression ``covariates $u, \ v$ match''  means  that the units with these covariates match.

 
{\it Assumption 1} (Stable Unit Treatment Value Assumption ({\it SUTVA}), Rubin, 1980, 1990) 
The distribution of potential outcomes for one unit is assumed to be independent of potential treatment status of another unit given the
observed covariates.

{\it Assumption 2} (Strong ignorability of treatment assignment given $x,$  Rosenbaum and Rubin, 1983)\\ 
{\it (i)} ${\cal R}$ and $T$ are conditionally independent given $x:$  ${\cal R} \ci  T \ |x,$ and\\
{\it (ii)} for every $t \in {\cal T},$ $0<p(t|x)$ (or equivalently $0<p(x|t)$).

Recall that 
$b(x)$ is a balancing score if  the conditional distribution of $x$ given $b(x)$ is the same for all treatment values, i.e.
\begin{equation}
 p(x|t,b(x))=p(x|b(x)), \mbox{ for all } t \in {\cal T}.
\label{eq:suf1}
\end{equation}
From (\ref{eq:suf1}), thinking of $t$ as  parameter value for the distribution of $x$ 
it follows that $b(x)$ is a sufficient statistic for the family ${\cal D}_{\cal T}=\{p(x|t); t \in {\cal T}\}.$  


\section{Matching and  Causal Inference with  $s$}

\quad 
 In this section 
the minimal sufficient statistic $s(x)$ is assumed  known. This is possible for various models.
 The results  are also applicable  for large samples when $s$ is estimated.

\subsection{Vector valued $s$}

\quad The first  result,  obtained directly from statistical theory, extends   Nelson and Noorbaloochi (2009, p. 619, Theorem 1) 
and justifies the use of $s$ for countably finite or countably infinite treatments. It also indicates that for multiple treatments 
$s$ will be often vector valued. 

\bep
\label{p:mss}
Assume that  ${\cal T}$ consists of either countably finite or countably infinite treatments and that the covariates' distributions
${\cal D}_{{\cal T}}$
 have all common support. Let ${\cal T}_k=\{t_1,\ldots, t_k\} \subset {\cal T}$ and let
\begin{equation}
\label{eq:mss1}
s(x)=s^{(1)}(x)\footnote{In $s^{(1)}(x), \ (1)$ indicates the denominator 
is $p(x|t_1).$}
=\left ( \frac{p(x|t_2)}{p(x|t_1)},  \frac{p(x|t_3)}{p(x|t_1)},  \ldots,  \frac{p(x|t_k)}{p(x|t_1)}\right ) 
\end{equation}
When $s$ 
is sufficient for ${\cal T}$ (i.e.  for ${\cal D}_{\cal T}$) it is also minimal sufficient. 
 \enp

\ber
\label{r:xvaluesdimension}
{\em  When 
${\cal D}_{{\cal T}_k}=\{p(x|t),  t \in {\cal T}_k\}$ do not have common support,  
$s(x)$'s 
 dimension 
depends on the $x$-values  (Lehmann and Casella, 1998, p. 70, Theorem 9.1). }
\enr  


{\bf The $s$-Matching Rule for covariates:}
 {\it Match $u$ to $v$ when $s(u)=s(v).$}

The $s$-matching rule can be used for  any number of units using their covariates.

\quad The next proposition shows that $s$-matching is not changed
when
\begin{equation}
\label{eq:mssj}
s^{(j)}(x)=\left ( \frac{p(x|t_1)}{p(x|t_j)},\ldots,  \frac{p
(x|t_{j-1})}{p(x|t_j)},
\frac{p(x|t_{j+1})}{p(x|t_j)}, \ldots,  \frac{p(x|t_k)}{p(x|t_j)}\right ), \\ j
\neq 1, 
\end{equation}
is used instead of $s=s^{(1)}$ in (\ref{eq:mss1}).


\bep
\label{p:equivs}
If $s(u)=s(v),$ then $s^{(j)}(u)=s^{(j)}(v), \ j>1.$
\enp
 Without loss of generality $s(x)=s^{(1)}(x)$ is used in this section.
 Propositions \ref{p:mss} and Remark \ref{r:xvaluesdimension}
 indicate clearly that with several treatments  $s$  is  not  expected to be scalar.
Nelson and Noorbaloochi (2009, p. 619) point out that $s(x)$ may be infinite dimensional thus contradicting
the existence of a finite dimensional parameter $\theta$ in {\em Assumption 3}
(Imai and Van Dyk, 2004). In Noorbaloochi,
Nelson and Asgharian (2010, p. 8, lines 1, 2) it is also mentioned that  ``{\em In general,  there is no univariate propensity score}.''   
This is indirectly confirmed for several data sets with high dimensional $x$-covariates for which the scalar propensity score deteriorates more
as $x$'s dimension increases (King et al., 2011, p. 18).


We revisit an example in Rosenbaum and Rubin (1983, p. 47) when the number of treatments $k$ is larger than 2.

 
\beex
{\em Let $p(x|t)$ be a polynomial exponential family distribution,
$$p(x|t)=h(x) \exp \{P_t(x)\}, \ t=1,\ldots, k,$$
with $P_t(x)$ a degree $m$ polynomial. Then, the statistic
\[
\left ( 
\ln \frac{p(x|t_2)}{p(x|t_1)}, \ldots, \ln \frac{p(x|t_k)}{p(x|t_1)}
\right)
=(P_2(x)-P_1(x), \ldots, P_k(x)-P_1(x)) 
\]
$$=\left (Q_1(x),\ldots,Q_{k-1}(x) \right)$$
is equivalent to the minimal sufficient statistic (\ref{eq:mss1}) with 
$Q_i(x)$ a degree $m$ polynomial, $i=1,\ldots, k-1.$}
\enex

\subsection{Causal inference framework and $s$}

\quad To obtain $s$ using the likelihood ratios in  (\ref{eq:mss1}) 
 the densities of the covariates
in  ${\cal D}_{\cal T}=\{p(x|t), \ t \in {\cal T}\}$ 
 have 
either to be  known or to be estimated. This may not be possible in practice.
Results are now presented to determine $s$ with the causal inference framework and without using ${\cal D}_{\cal T},$ simply  from the conditional density $q(t|x)$  of $T$ given the $x$-covariates.
The first result  involves  ${\cal D}_{\cal T}$ but it  is used to prove subsequent results.

\bep (see, e.g. Chen, 2010, Ch. 6)  Let ${\cal D}_{\cal T}$ 
be the family of the $x$-covariates 
densities. Assume that there exist function $s^*(x)$ such that for any covariates $x_1$ and $x_2$ 
 the ratio $\frac{p(x_1|t)}{p(x_2|t)}$ is constant as function of $t$  if and only if $s^*(x_1)=s^*(x_2).$  Then, $s^*(x)$ is 
minimal sufficient. 
\enp


The tool to determine $s$  via $q(t|x)$  is  the decomposition
\begin{equation}
\label{eq:Bmsf}
p(x|t)=q(t|x) \cdot p_X(x) \cdot \pi^{-1}_T(t),
\end{equation}
that leads to $T$'s posterior factorization criterion (PFC) and the coarsest balancing score criterion (CBSC).
\bep
\label{p:Bmsf}
Let $X$ and $T$ be random vectors in Euclidean spaces with densities, respectively, $p_X$ and  $p_T$ and with conditional densities  
$p(x|t)$ and $q(t|x).$ Then,\\
a) (Posterior Factorization Criterion) $s(x)$ is sufficient statistic for ${\cal D}_{\cal T}$  if and only if 
\begin{equation}
\label{eq:BNfact}
q(t|x)=g_1(s(x),t)\cdot  g_2(t) \cdot g_3(x) \ \forall \ x, \ t.
\end{equation}
b) (Coarsest Balancing Score Criterion) Assume that for any $x_1$ and $x_2$  the ratio $\frac{q(t|x_1)}{q(t|x_2)}$ is independent of $t$ 
if and only if $s^*(x_1)=s^*(x_2).$ Then, $s^*$ is minimal sufficient statistic for ${\cal D}_{\cal T}.$
\enp

 Proposition \ref {p:Bmsf} is  used to derive directly  previous results in the literature.

\beex (The  propensity score, Rosenbaum and Rubin, 1983) {\em The treatments ${\cal T}=\{1,2\}$ and the propensity score 
$e(x)=q(1|x).$ For $q(t|x)$ it holds 
\begin{equation}
\label{eq:qt}
q(1|x)=e(x), \ q(2|x)=1-e(x).
\end{equation}
From (\ref{eq:qt}) and Proposition  \ref {p:Bmsf} $a)$ 
$s(x)=e(x)$ is sufficient statistic.
Since 
the ratio 
$\frac{q(t|x_1)}{q(t|x_2)}$ is independent of $t$ for all $x_1, \ x_2$  if and only if
$$\frac{q(1|x_1)}{q(1|x_2)}=\frac{q(2|x_1)}{q(2|x_2)} \hspace{3ex}  \Leftrightarrow  \hspace{3ex} e(x_1)=e(x_2),$$
 from Proposition \ref {p:Bmsf} $b)$
$e(x)$ is minimal sufficient. The same result is obtained via 
Proposition \ref{p:mss}.}
\enex

\beex (The propensity function, Imai and van Dyk, 2004) {\em For various kinds of treatments $t \in {\cal T},$  the propensity function $e_{\psi}(\cdot |x)=q_{\psi}(\cdot |x)$ depends on $x$ only through the unique,  finite dimensional parameter $\theta_{\psi}(x) (Assumption \ 3); \psi$  is known parameter. Therefore,
$q_{\psi}(t|x)$  has form (\ref{eq:BNfact}) with $g_3(x)=1$ for every $x$ and from Proposition  \ref {p:Bmsf} $a)$  $s(x)=\theta_{\psi}(x)$ is sufficient.}
\enex
 
 Proposition \ref {p:Bmsf}  extends the results in Imai and Van Dyk (2004) by weakening their {\em Assumption 3} 
to accommodate more $q(t|x)$-models.

\bec 
\label{c:extImVD} Let $q_{\psi}(t|x)$ be the density of $T$ given the $x$-covariates;  $\psi$ are known parameters. Assume that  there are functions $\theta_{1, \psi}(x),  \ g_{1, \psi}, \ g_{2, \psi}, \ g_{3,\psi}$ such that 
\begin{equation}
\label{eq:extImVD}
q_{\psi}(t|x)=g_{1,\psi}(\theta_{1,\psi}(x), t)\cdot g_{2,\psi}(t) \cdot g_{3,\psi}(x) \ \forall \ x, \ t.
\end{equation}
Then,\\
a) $\theta_{1,\psi}(x)$ is a balancing score, and\\
b) $\theta_{1,\psi}(x)$ is the coarsest balancing score when for every $x_1, \ x_2,$ the ratio
$$\frac{g_{1,\psi}(\theta_{1,\psi}(x_1), t)}{g_{1,\psi}(\theta_{1,\psi}(x_2), t)}$$  is independent of $t$ if and only if
$\theta_{1,\psi}(x_1)=\theta_{1,\psi}(x_2).$
\enc

The minimal sufficient statistic $s$ for ${\cal D}{\cal T}$ is now determined when $T$'s posterior is a generalized linear model.

\beex {\em Assume the treatment variable $T$ with values in $R^d$  is modeled given the $x$-covariates ($\in R^p$) with  a generalized linear model in canonical form, i.e.
\begin{equation}
\label{eq:glm1} 
q(t|x)=\exp\{t' b(x)+c(x)+d(t)\};
\end{equation} 
$t'$ denotes $t$'s transpose and $n \in R^d.$
From Proposition   \ref {p:Bmsf} a),  $b(x)$ is sufficient statistic for ${\cal D}_{\cal T}.$ 
Since the ratio 
$$\frac{q(t|x_1)}{q(t|x_2)}=exp\{t'[b(x_1)-b(x_2)]+c(x_1)-c(x_2)\}$$
is independent of $t$ for all $t, x_1,x_2$  if and only if $b(x_1)=b(x_2),$
from  Proposition \ref {p:Bmsf} b)  $b(x)$ is also minimal sufficient.

There are different forms $b(x)$ can have. For example, when $b$ takes real values, 
\begin{equation}
\label{eq:glmmss}
b(x)=\phi(\sum_{j=0}^p \beta_jx_j), \hspace {4ex} 
b(x)=\sum_{j=0}^K\beta_j  \phi_j(x);
\end{equation}
$x_0=1, \ x_j$ is $x$'s $j$-th coordinate, $\beta_j \in R,$
the functions $\phi, \ \phi_j$
are  assumed to be known and real valued, $j$'s values are according to the corresponding sum.

When
$b(x)$ is  known 
it can be used for matching units from different populations.
When $b(x)$ is not known, it has to be estimated with $\hat b(x)$ that is used for matching.}
\enex

\ber
\label{r:glmviol}
{\em For the generalized linear model (\ref{eq:glm1}) {\em Assumption 
3} in Imai and Van Dyk (2004) does not hold because $q(t|x)$ depends  on $x$ via $c(x)$ also.}
\enr

\subsection{Simultaneous causal comparisons}

\quad The key result allowing for simultaneous causal comparison of several treatments follows,
establishing  strong ignorability of treatment assigment given $s(x).$ 


\bep 
\label{p:sita}
Under Assumption 2, for the responses ${\cal R}$ 
and the treatment variable $T=t$ it holds 
\begin{equation}
\label{eq:sita}
p\{t, {\cal R} | s(x)=s\}=p\{t|s(x)=s\} \cdot p\{ {\cal R}| s(x)=s\}.
\end{equation}
\enp

Proposition \ref {p:sita} suggests simultaneous causal comparisons 
using $s(x)$  to balance subpopulations  for all treatments
and obtain unbiased estimates of the average treatment effects.

\bep
\label{p:unbs}
 Suppose that treatment assignment is strongly ignorable (Assumption 2) and that a value $s_0$ of $s(x)$ is randomly sampled from the population of units with
covariates $x \in {\cal C}.$ 
Units receiving treatments $t_i$ and $t_j$ are sampled with 
$s$-value for their covariates equal $s_0, \ i \neq j.$ 
Then, the expected difference in response for 
the units chosen is the expected treatment effect at $s(x)=s_0.$ 
The mean of such pair differences 
over all $s(x)$-values is unbiased for the average treatment effect 
$E\{r(t_i)-r(t_j)\}$ and the same holds, concurrently  given $s(x),$ for any number of average treatment effects. 
\enp


\section{Implementation}

\quad For the
$s$-matching's implementation  there are practical issues some of which depend 
on the data to be analyzed and the
 assumptions on the $x$-covariates models. Among these issues {\em a)} $s$ should be determined, {\em b)}  when  likelihood ratios have to be estimated,  the curse of dimensionality  problem  should be addressed,  and  {\em c)}  the dimensionality of $s$ may  be reduced if there is no much loss of information. Some directions for {\em a)}-{\em c)} follow.

 Known theorems in statistics  (see, e.g.,   Lehmann and Casella, 1998) 
 allow to obtain the minimal sufficient statistic  $s.$ 
 Proposition 3.1 is used with a small number  $k$ of treatments  to define a minimal sufficient statistic  $s^*$ for the corresponding distributions, ${\cal D}_{{\cal T}_k},$ and then show that $s^*$ is sufficient for all the
 $x$-covariates distributions,  ${\cal D}_{{\cal T}},$ i.e. $s=s^*.$  For example,   if $p(x|t)$ follows a normal distribution with mean $t$ and known variance
(say) 1, 
the  minimal sufficient statistic, $s^*,$ is determined for the distributions ${\cal D}_{\{t_1,t_2\}}; \ t_1, \ t_2$ are treatments, $t_1 \neq t_2.$
 Neyman's factorization criterion shows that $s^*$ is sufficient for all
 $t$-values, so it is minimal sufficient for  ${\cal D}_{{\cal T}}.$  Alternatively, when $p(x|t)$ belongs to a $p$-parameter exponential family in canonical form,  Neyman's factorization criterion determines  $s(x)$ that is also minimal sufficient if the parameter space ${\cal T}$ contains an open, $p$-dimensional rectangle.


 For the implementation of Proposition 3.1 in applications  only a subpopulation $\tilde {\cal P}_t$ of  ${\cal P}_t$ may be  available.  Let ${\cal MP}_t$ denote the units to be matched from 
$\tilde {\cal P}_t$-subpopulation. Use
$s=s^{(1)}$ in (\ref{eq:mss1}) 
 to
match 
a unit in ${\cal MP}_t$
having covariates $u$
with a unit from $\tilde {\cal P}_r$
having 
covariates $v_{m,r} \in {\cal C}(\tilde {\cal P}_r),$ such that
\begin{equation}
\label{eq:practicalcriterion}
v_{m,r}=\arg \min_{v \in {\cal C}(\tilde {\cal P}_r)} ||s(u)-s(v)||^2, \ r \in T-\{t\};
\end{equation}
$|| \cdot ||$ is the usual Euclidean distance in $R^p$ and   in $v_{m,r}$ the index  $m$ denotes ``matching'' unit  from sub-population 
$\tilde {\cal P}_r.$ This approach  is the nearest neighbor $1:1$  matching with replacement and can be properly modified for  
$1:k$  matching with or 
without replacement. For more information on  matching methods and for optimal matching questions  see, e.g., Rosenbaum (1989) and  Stuart (2010).

Additional matching sets 
for ${\cal MP}_t$ can  be  obtained using $s=s^{(j)}$ 
(or its estimates)  in (\ref{eq:practicalcriterion}), $j=2,\ldots,k,$ and
the decision maker can select the ``best'' matching set, for example, that with
the nearest means to the ${\cal MP}_t$ covariates' means 
with respect to $|| \cdot ||$ or
the sup-norm distance $|| \cdot ||_{\infty}.$

When the form of the $x$-covariates densities in ${\cal D}_{\cal T}$ is not known and ${\cal T}={{\cal {T}}_k},$ the usual nonparametric estimation of each density in the ratios (\ref{eq:mss1})  is affected by the curse of dimensionality of the $x$'s. Rather
than estimating separately  each density one may use the approach adopted  in Machine Learning
 for determining the ratio of the densities from the training and test data. The densities' ratio  is  expressed as linear model 
with respect to a basis of functions. The coefficients are estimated using  observations  from the two populations according to a given method that usually ends with a  convex minimization problem. For the description of the 
estimation methods  in Machine Learning see, e.g., Sugiyama {\em et al.} (2007) and Nguyen {\em et al.} (2010).
 For independent $x$-samples from nonparametric models, one  for each $t \in {\cal T},$ the class of empirical distributions
is minimal sufficient. 

 The use   of $s$  is supported by the  existing statistical theory that allows 
{\em i)} to obtain approximate sufficient statistics  and evaluate the approximation's error (see, e.g., Le Cam, 1964, 
Joyce  and Marjoram, 2008), and
 {\em ii)}   use a principal components transformation of its components to investigate whether a sufficient summary of smaller
dimension exists  (Nelson and Noorbaloochi, 2009). 

\section{Appendix}

\quad {\bf Proof of Proposition \ref{p:mss}:}  It is direct consequence of Theorem 6.12, in Lehmann and Casella, 1998, p. 37
and related theorems therein.\\

{\bf Proof of Proposition \ref{p:equivs}:} Since $s(u)=s(v),$ it holds
 \begin{equation}
\label{eq:equal1}
 \frac{p(u|t_i)}{p(u|t_1)}=\frac{p(v|t_i)}{p(v|t_1)}, \ i=2,\ldots,k.
\end{equation}
In (\ref{eq:equal1}), divide  the $i$-th equality with the $j$-th equality, $i\neq j ,$ and invert the $j$-th equality to obtain
$$ \frac{p(u|t_i)}{p(u|t_j)}=\frac{p(v|t_i)}{p(v|t_j)}, \ i \neq j, \mbox{ or } s^{(j)}(u)=s^{(j)}(v).$$



{\bf Proof of Proposition \ref {p:Bmsf}:}  {\em a)} From Neyman's Factorization criterion $s$ is sufficient statistic if and only if
\beeq
\label{eq:NFC}
p(x|t)=h_1(s(x),t) \cdot h_2(x) \ \forall \ x, \ t,
\eneq
and from  decomposition (\ref{eq:Bmsf}) 
$$\Leftrightarrow  \hspace{3ex}
q(t|x)\cdot p_X(x)\cdot \pi^{-1}_T(t)=h_1(s(x),t) \cdot h_2(x) 
$$
$$
 \Leftrightarrow \hspace{3ex}  q(t|x)= h_1(s(x),t) \cdot \pi_T(t)  \cdot h_2(x)\cdot p^{-1}_X(x).$$
Equality  (\ref{eq:BNfact}) follows with
$$g_1(s(x),t)=h_1(s(x),t), \hspace{4ex}   g_2(t)=\pi_T(t),  \hspace{4ex}  g_3(x)= h_2(x) \cdot p^{-1}_X(x).$$
Conversely, from (\ref{eq:BNfact}) Neyman's Factorization criterion  (\ref{eq:NFC}) is obtained via (\ref{eq:Bmsf}).

{\em b)} From decomposition (\ref{eq:Bmsf})  the ratio 
$$\frac{p(x_1|t)}{p(x_2|t)}=\frac{q(t|x_1)}{q(t|x_2)} \cdot  \frac{p_X(x_1)}{p_X(x_2)}$$
is independent from $t$ if and only if the ratio
$$\frac{q(t|x_1)}{q(t|x_2)}$$
is independent of $t$ and this holds if and only if 
$$s^*(x_1)=s^*(x_2).$$
Thus, from Proposition 4.1 $s^*$ is minimal sufficient statistic. 

{\bf Proof of Corollary  \ref{c:extImVD}:} Both parts follow from Proposition \ref {p:Bmsf}.

{\bf Proof of Proposition \ref{p:sita}:} The proof follows the lines in Imai and Van Dyk (2004),
$$p\{x, t, {\cal R}|s(x)=s\}=p\{x,t|s(x)=s\}\cdot p\{{\cal R}|x, t, s(x)=s\}$$
$$=p\{t|s(x)=s\}\cdot  p\{x|t, s(x)=s\}  \cdot p\{{\cal R}|x, t, s(x)=s\}$$
$$=p\{t|s(x)=s\}\cdot  p\{x| s(x)=s\}\cdot p\{{\cal R}|x, s(x)=s\}.$$
The third equality is obtained using  Proposition \ref{p:mss} and strong ignorability of treatment assignment given $x$
(Assumption 2). It follows that
$$p\{t, x, {\cal R}|s(x)=s\}=p\{t|s(x)=s\}\cdot  p\{x, {\cal R}|s(x)=s\}$$
Integrating  both sides of the last equation over the $x$'s for which $s(x)=s,$
we obtain that given $s(x)=s,$ ${\cal R}$ and $T$ are independent.

{\bf Proof of Proposition \ref{p:unbs}:} From Assumption 2,
$$E \{r(t_i)|s(x)=s, T=t_i  \}-E\{r(t_j)|s(x)=s, T=t_j\}$$
$$= E\{r(t_i)|s(x)=s\}-E\{r(t_j)|s(x)=s\}=E\{r(t_i)-r(t_j)|s(x)=s\}$$
and it follows that
$$E_{s}\left[ E\{r(t_i)-r(t_j)|s(x)=s\}\right]=E\{r(t_i)-r(t_j)\};
$$
$E_s$ denotes expectation with respect to all values $s$ of $s(x), \ x \in {\cal C}.$


\end{document}